\documentclass[10pt, runningheads]{llncs}
\usepackage{amsmath}
\usepackage{amssymb}
\usepackage{color}
\usepackage{fancyhdr}
\usepackage{graphicx}
\newcommand{\pv}{\par\vspace{1ex}}

\newcommand{\N}{{\mathbb{N}}}
\newcommand{\R}{{\mathbb{R}}}
\newcommand{\Q}{{\mathbb{Q}}}

\newcommand{\F}{{\mathbb{F}}}

\begin{document}

\title{On diagrams accompanying \textit{reductio ad absurdum} proofs in Euclid's \textit{Elements} 
book I. {Reviewing Hartshorne and Manders}  }

\author{Piotr Błaszczyk\inst{1}\orcidID{0000-0002-3501-3480} \and \\ Anna Petiurenko\inst{1}\orcidID{0000-0002-0196-6275} %\and
%Third Author\inst{3}\orcidID{2222--3333-4444-5555}
}
\authorrunning{Piotr Błaszczyk, Anna Petiurenko}
\institute{Institute of Mathematics, Pedagogical University of Cracow, Podchorazych 2, Cracow, Poland\\
\email{piotr.blaszczyk@up.krakow.pl;\ anna.petiurenko@up.krakow.pl}\\
\url{https://matematyka.up.krakow.pl/} 
%\and
}  %Podchorazych 2, Cracow,
\maketitle              % typeset the header of the contribution
\begin{abstract} 
 Exploring selected  \textit{reductio ad absurdum} proofs in Book 1 of the \textit{Elements}, we show they include figures that are not constructed. It is squarely at odds with Hartshorne's claim that ``in Euclid's geometry, only those geometrical figures exist that can be constructed with ruler and compass".
 
We also present diagrams questioning Manders' distinction between exact and co-exact attributes of a diagram, specifically,  a model of semi-Euclidean geometry which satisfies   \textit{straightness of lines} and \textit{equality of angles} and does not satisfy the parallel postulate.

%Euclidean diagrams  \and \and Cartesian plane over Euclidean field
\end{abstract}

\section{Introduction}
Euclid's propositions are of two kinds: constructions and demonstrations.
 I.1 and I.32 are model examples: the first requires the construction of an equilateral triangle, the second --  a demonstration that angles in a triangle sum up to {``two right angles" ($\pi$, in short)}. Yet,  the proof of I.32  includes the construction of a parallel through a point. 
Indeed,  each proposition includes a construction part (\emph{kataskeu\={e}})  which introduces auxiliary lines exploited in the proof (\emph{apodeixis}). In this paper, we focus on diagrams accompanying \textit{reductio ad absurdum} proofs, as they undermine the common understanding of Euclid's diagram.

%Indeed, the text of each proposition is made up of  six parts: \emph{protasis} (stating the relations among geometrical objects by means of abstract and technical terms), \emph{ekthesis} (identifying objects of \textit{protasis} with lettered objects), \emph{diorisomos} (reformulating \textit{protasis} in terms of lettered objects), \emph{kataskeu\={e}}   (a construction part which introduces auxiliary lines exploited in the proof that follows),  \emph{apodeixis} (proof, which usually proves  the \emph{diorisomos}' claim),  \emph{sumperasma} (reiterating \textit{diorisomos}).  
%In this sense, every proposition includes a construction, meaning \emph{kataskeu\={e}} part.
%In this paper, we focus on diagrams accompanying \textit{reductio ad absurdum} proofs, for they undermine the common understanding of Euclid's diagram. % is the most insightful reading of the \textit{Elements} ever. It

{The 20th-century} Euclid scholarship grows in two trends: mathematical and historical-philological. For diagrams, 
 the first school examines them only as straightedge and compass constructions, while the second seeks to show they convey some mathematical information beyond construction requirements. We challenge both approaches with specific diagrams.

R. Hartshorne \cite{ref_RH} 
 develops a coherent reading of the {\textit{Elements}, Books I--IV} focused on tacit axioms, non-defined concepts, or relations and interprets them in the system of Hilbert axioms.
He does not find Euclidean diagrams problematic, misleading, or competing with a logical account of geometry;  on constructions, though, he writes: 
``The constructive approach pervades Euclid's \textit{Elements}. There is no figure in the entire work that cannot be constructed with ruler and compass 
 [...] in Euclid's geometry only those geometrical figures exists that can be constructed with ruler and compass"
 (\cite{ref_RH}, 18--19). 
 
We discuss two examples undermining Hartshorne's claim: the figure accompanying proposition I.7 and one implied by the proof of I.27.  The first is not, and indeed, cannot be constructed, as assumptions of the proposition introduce an inconsistent object. 
The non-constructive  mode of the second figure is related to the requirement ``being produced to infinity" inherent in the definition of parallel lines.  
These two figures are by no means incidental, as the first props the SSS theorem (I.8),  and the second brings us to the core of the Euclid system.

J. Ferreir{\'o}s  provides a concise picture of the second school: ``The original geometry of Euclid
lacks the means to derive its theorems by pure logic, but it presents us with a most interesting and fruitful way of proving results by diagrams" (\cite{ref_JF}, 132).\ 
%Unlike Hartshorne,  the second school is so far at the stage of justifying selected propositions rather than the %Euclid system as a whole. 
In this vein, K. Manders \cite{ref_KM} introduced a distinction between \textit{exact} and \textit{co-exact} information inferred from a diagram that got widespread renown among scholars exercising diagrammatic approach. 
Respective definitions read:        
``Co-exact attributes
are those conditions which are unaffected by some range of every continuous
variation of a specified diagram; paradigmatically, that one region includes
another [...], %(which is unaffected no matter how the boundaries are to some extent shifted and deformed)
 or the existence of intersection points such as
those required in Euclid I.1 [...]. 
%(which is unaffected no matter how the circles are to some extent deformed). 
Exact attributes are those which, for at least
some continuous variation of the diagram, obtain only in isolated cases;
paradigmatically, straightness of lines or equality of angles [...]. Exact attributes [...] are unstable under perturbation of a
diagram  (\cite{ref_KM}, 92--93).\footnote{\cite{ref_KM}, p. 92 presents parallelism as an exact attribute.}
%Specifically, one can read off from a diagram that a region includes another, or lines meet, while should not infer equality of angles or that lines are parallel.

{Instead of `continuous variation of a diagram', we introduce a global perspective, meaning specifications of the plane on which a diagram lies. Accordingly, we examine  Euclid diagram I.1  on various Cartesian planes showing that the existence of the intersection of circles involved depend on  characteristics of a plane. 
Regarding \textit{exact} attributes, we present a model of 
a semi-Euclidean plane that does not affect \textit{straightness of lines or equality of angles} but affects parallelism (especially I.29). 
Both counterexamples meet the scheme: without touching a diagram but changing assumptions on the space hosting it, we get different results concerning \textit{co-exact} (intersection of circles) and \textit{exact} (parallelism) attributes. 
Depending on assumptions concerning space, the same (from the diagrammatic perspective) circles meet or not, and the same straight lines are parallel or not.

 As for co-exact attributes, we  provide an analysis of proposition I.6 that undermines Mandres' interpretation of inequality in terms of part-whole.}

For the most part, our arguments exploit an interpretation of \textit{greater-than} relation. It affects  our account of the deductive structure of the \textit{Elements},  some existential claims,  Manders' interpretation of \textit{greater-than} in terms of part-whole, and his claim regarding exact attributes. 
 {Euclid's arguments exploring that
 relation proceed \textit{reductio ad absurdum} mode. Significantly, out of eleven indirect proofs in Book I, ten, specifically I.6, 7, 27, 29, employ \textit{greater-than} relation}.\footnote{{The others are I.14, 19, 25, 26, 39, 40.}}

% \footnote{Cf. ``Bringing Euclid up to Hilbert’s standard means banishing diagrams from the proofs and replacing them with an abstract theory of order. The evidence for any such theory in the \textit{Elements} is nil" \cite{ref_JM}, 274.}

\section{Euclid and Hilbert construction tools}

Hartshorne (\cite{ref_RH}, 102) introduces term \textit{Hilbert construction tools}, meaning  transportation of line segments and angles. Hilbert axioms justify these operations; besides, they also state the uniqueness of respective line segments and angles.\footnote{\cite{ref_MG}, pp. 597--602 provide a concise account of Hilbert axioms. C1 and C4, respectively, decree construction tools.} 
Diagrams drawn up with both tools are acquired using the first alone; it suggests Euclid's straightedge and compass are more effective. We follow that clue to contrast Euclid and Hilbert approaches.

\subsection{Transportation of segments. I.1--3}
I.1  \textit{To construct an equilateral triangle on the given line AB}.\footnote{{English translation by Fitzpatrick \cite{ref_fitz}, diagrams after \cite{ref_JLH}.}}
\begin{figure}[h!]
\begin{center}
\includegraphics[scale=0.7]{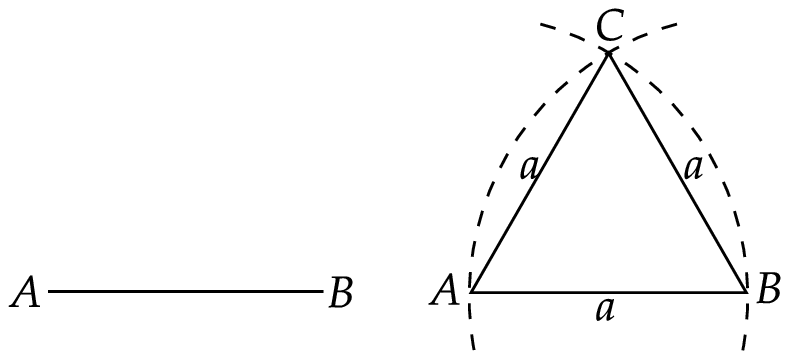}
%\caption{Proof of the \textit{Elements}, I.1 schematized} 
\label{figI1}
\end{center}
\end{figure}
 
 %\centerline{\includegraphics[scale=0.7]{I_1a}}
Given that $AB=a$, point \textit{C}, the third vertex of the wanted triangle is an intersection of circles $(A, a)$ and $(B, a)$.    
 In tables like the one below, we lay out points resulting from intersections of straight lines and circles.\footnote{The idea of such tables originates from \cite{ref_GM}.} 
\pv
 \begin{table}[!h]
 \centering
 \begin{tabular}{c }
 $(A, a), (B, a)$ \\ \hline
 $C$  \\  
 \end{tabular}
 \end{table}

I.2 \textit{To place a straight-line at point A 
equal to the given straight-line BC} [$b$].
\begin{figure}[h!]
\begin{center}
\includegraphics[scale=0.7]{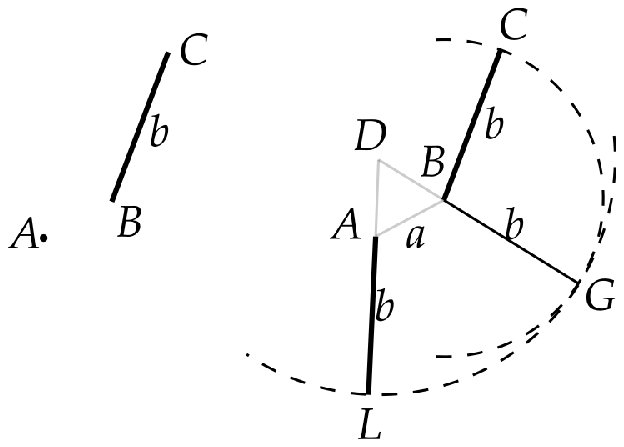}
%\caption{Proof of the \textit{Elements}, I.2 schematized} 
\label{figI2}
\end{center}
\end{figure}

On the line-segment \textit{AB}, we construct an equilateral triangle \textit{ABD} with side \textit{a}; the below diagram depicts {its \textit{shadow} in grey}. Point \textit{G} is the intersection of the circle $(B,b)$ and the half-line 
  $DB^{\rightarrow}$. Now, \textit{DG} represents the sum of line-segments $a,\,b$.
  %\footnote{Clearly,   all throughout the \textit{Elements} there are no term-counterpart of the modern concept \textit{sum}, or \textit{addition}.} 
  Circle $(D,a+b)$ intersects the half-line 
  $DA^{\rightarrow}$ at point \textit{L}. Due to the \textit{Common Notions} (CN, in short) 3, \textit{AL} proves to be equal $b$. 
   
  \pv
\begin{table}[!h]
 \centering
 \begin{tabular}{c|c}
$(B, b), \,\, {DB}^{\rightarrow}$  & $(D, a+b), \,\, {DA}^{\rightarrow}$ \\ \hline
 $G$				    & $L$\\
 \end{tabular}
 \end{table}

 %\centerline{\includegraphics[scale=0.7]{I_2a}}

 Owning to I.1-2, $b$ is placed at \textit{A} in a precise position. Drawing  circle \textit{(A,b)}, one can choose any other position at will,  and that is the substance of I.3.
 %what Euclid does in the subsequent proposition.
  \pv 
  I.3 \textit{To cut off  a straight-line equal to the lesser C} [$b$] \textit{from the greater AB.}
  
  \begin{figure}[h!]
\begin{center}
\includegraphics[scale=0.7]{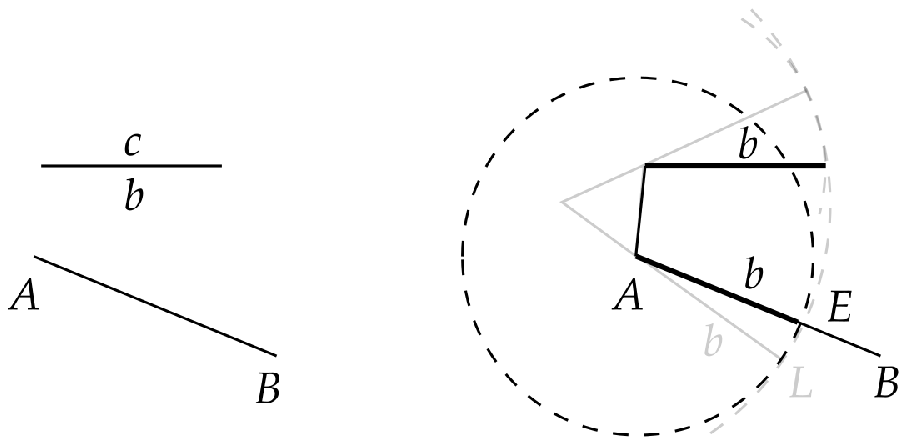}
%\caption{Proof of the \textit{Elements}, I.3 schematized} 
\label{figI3}
\end{center}
\end{figure}
  
    %\centerline{\includegraphics[scale=0.7]{I_3b}}
       
   Line-segment $b$ is transported to \textit{A} into position \textit{AL}; the above diagram depicts the \textit{shadow} of that construction; let $Ab$ be its symbolic representation. The intersection of circle $(A, b)$ and line \textit{AB} determines \textit{E} such that ${AE=b}$.

   \begin{table}[!h]
 \centering
 \begin{tabular}{c|c }
 $Ab$ &  $(A, b),\,\, AB$ \\ \hline
 &$E$  \\  
 \end{tabular}
 \end{table}

Summing up, due to I.1--3, one can transport any line segment to any point and position. An equilateral triangle is a tool to this end, while the existence of circle-circle and circle-line intersection points are taken for granted.

The Euclid system requires a circle-circle or circle-line axiom, both finding grounds in \textit{Postulates} 1--3 that introduce straight-edge and compass. Logically, these two tools reduce to compass alone (\textit{vide} Mohr-Mascheroni theorem), yet, throughout the ages, the economy of diagrams prevailed and no one questioned the rationale for Euclid's instruments.
There are, however, models of the Hilbert system that do not satisfy the circle-circle axiom. Moreover, Hartshorne shows (\cite{ref_RH}, 373) that I.1 does not hold in the Hilbert system of absolute geometry.  
Thus, already at the very first proposition of the \textit{Elements}, we observe that Euclid and Hilbert's systems follow alternative deductive tracks. Therefore one cannot simply merge  Hilbert's axioms with Euclid's arguments.

\subsection{Transportation of angles. I.22--23}
In I.22, Euclid builds a triangle from three given line segments.\footnote{\cite{ref_MG}, p. 173 observes it is equivalent to the circle-circle axiom.}  

\centerline{\includegraphics[scale=0.7]{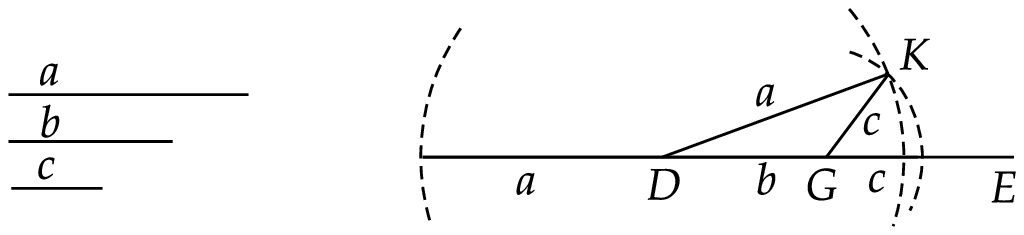}}

The below table presents  $D, E$ as random- and $G, K$ as intersection- points.

   \pv \begin{table}[!h]
 \centering
 \begin{tabular}{c|c|c|c|c|c}
              &${Db}$  & $(D, b),\, DE$& $Da$ & $Gc$ &$(D, a)$, $(G, c)$\\ \hline
 $D$, $E$&            &                $G$&         &   & $K$\\ 
 \end{tabular}
 \end{table}

  I.23, angle transportation, rests on triangle construction as follows: on sides of the given angle, Euclid builds a triangle, transports its sides to $A$, $G$, obtaining another triangle.  By the SSS, $\triangle CDE=\triangle AFG$, hence  $\angle KCL=\angle FAG$.%\footnote{Alternatively, one can  draw the perpendicular through $D$ and apply SAS.} 

\centerline{\includegraphics[scale=0.7]{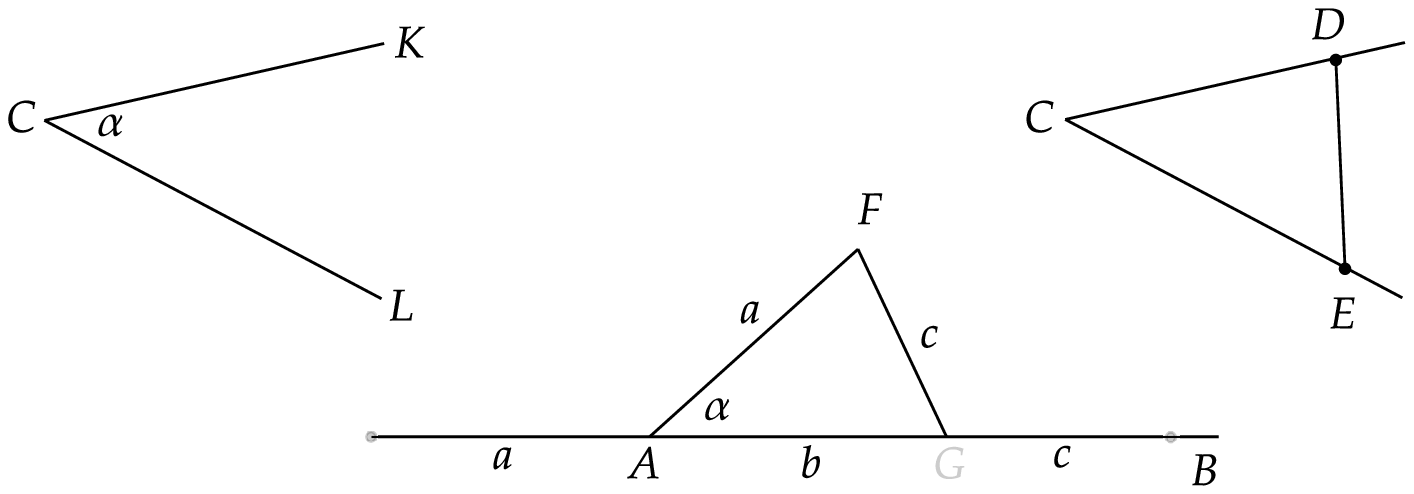}}

  %\centerline{\includegraphics[scale=0.7]{I_23}}
  %\centerline{\includegraphics[scale=0.7]{I_23a}}

    \begin{table}[!h]
 \centering
 \begin{tabular}{c|c|c|c|c}
 ${Ab}$  & $(A, b),\, AB$& $Aa$ & $Gc$ &$(A, a)$, $(G, c)$\\ \hline
             &                $G$&         &   & $F$\\ 
 \end{tabular}
 \end{table}
 
 \section{Euclid's vs Hilbert's deduction:   SAS to SSS}

Throughout propositions I.1--34, \textit{equality} means congruence, whether applied to line segments,  angles, or triangles. In I.5--8, showing the 
SSS theorem, Euclid assumes I.4,  \textit{Common Notions}, and characteristics of the \textit{greater-than} relation.

%\centerline{\includegraphics[scale=0.7]{I_4}}

The proof of I.4 (SAS criterion) relies on the \textit{ad hoc} rule: \textit{two straight-lines can not encompass an area}.
The diagram depicts an area encircled by the base  $EF$ of the triangle and a curve with ends $E, F$.  
By contrast,  Hilbert axioms guarantee a unique line through points $E,\, F$ and diagram I.4 has no grounds.

\pv I.5 \textit{Let ABC be an isosceles triangle. I say that the angle ABC is
 equal to ACB}.
 
  \centerline{\includegraphics[scale=0.7]{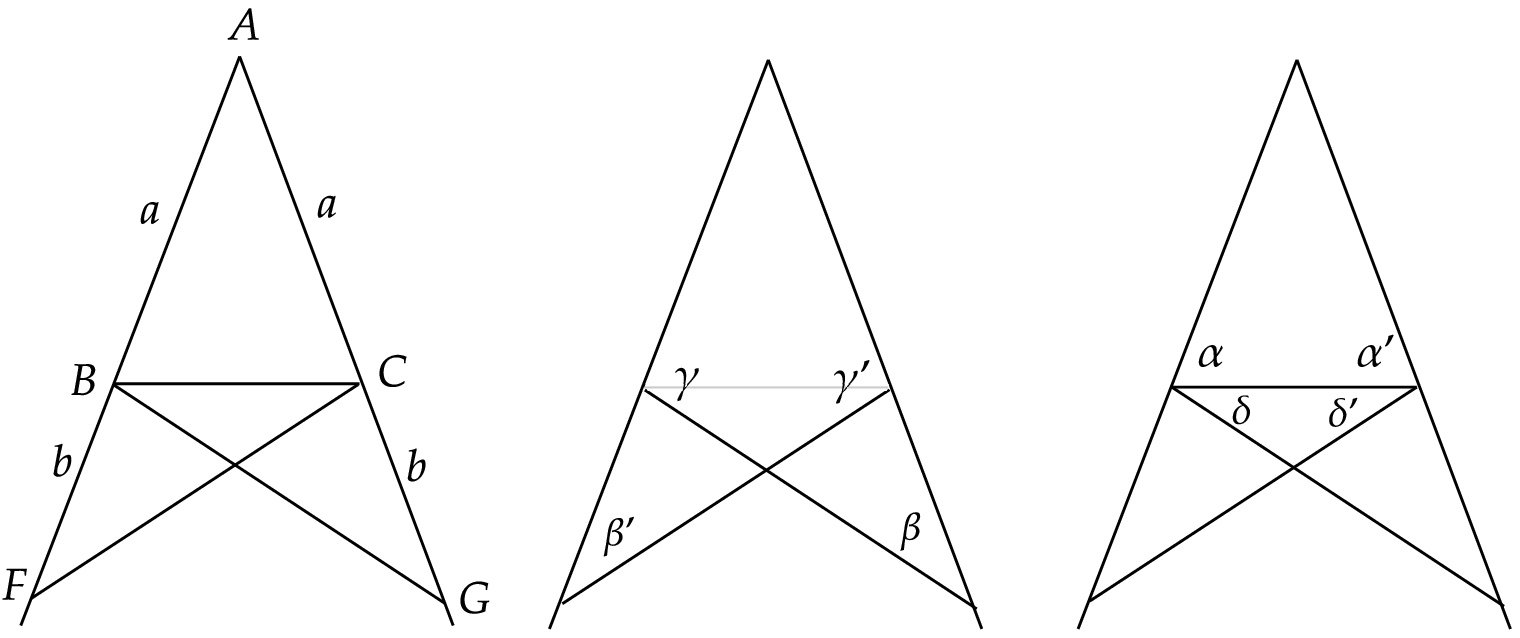}}
  
  The construction part is simple: \textit{F} is taken at random on the half-line $AB^{\rightarrow}$, then \textit{G} such that $AF=AG$ is determined on the half-line $AC^{\rightarrow}$. %The below table details the construction.  
  \pv
\begin{table}[!h]
 \centering
 \begin{tabular}{c|c}
$AB^{\rightarrow}$ &$(A, a+b), \,\, AC^{\rightarrow}$  \\ \hline
$ F$&$G$ \\ 
 \end{tabular}
 \end{table}
 
  (i) Now, due to SAS, $\triangle GAB=\triangle FAC$. Thus $FC=BG$ and
  $$\beta =\angle AGB=\angle AFC=\beta',$$ 
  $$\gamma=\angle ABG=\angle ACF=\gamma'.$$

 (ii) Again by SAS,  $\triangle BFC=\triangle BGC$, and 
 $$\delta=\angle CBG= \angle BCF=\delta'.$$ 
 
 %(iii)  $\angle ABC=\angle ABG-\angle CBG, \ \ \ \angle ACB=\angle ACF- \angle BCF.$
 
 (iii) By CN 3, $\gamma-\delta=\gamma'-\delta'$. Since
  \[\alpha=\gamma-\delta,\ \ \gamma'-\delta'=\alpha',\]
  
  the equality $\alpha=\alpha'$ holds.
 %\[\angle ABG-\angle CBG, \ \ \ \angle ACB=\angle ACF- \angle BCF.\]
 \hfill{$\Box$}

  \pv I.6 \textit{Let ABC be a triangle having the angle ABC equal
to the angle ACB. I say that side AB is also equal to side AC.}

\begin{figure}
\begin{center}
\includegraphics[scale=0.7]{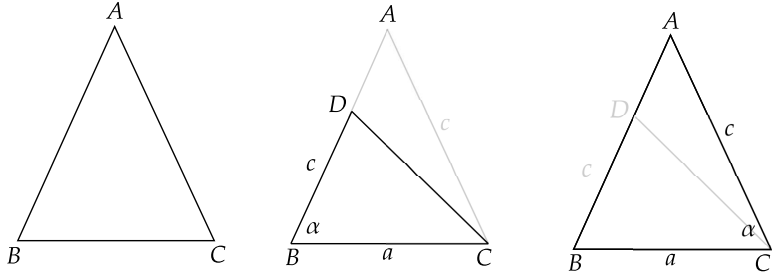}
\caption{Proof of I.6 schematized.} \label{figI6}
\end{center}
\end{figure}

 The proof reveals assumptions in no way conveyed through definitions or axioms. 
 At first, it is the trichotomy law for line segments. Let $AB=b$, $AC=c$ (see Fig. \ref{figI6}). To reach a contradiction Euclid takes: 
if $b\neq c$, then $b<c$ or $b>c$. Tacitly he assumes that exactly one of the conditions holds
\[b<c,\ \ \ b=c,\ \ \ b>c.\]

Let $b>c$. Then the construction follows: ``let DB, equal to the lesser AC, have been cut off from the greater AB". However, given that angles at $B$ and $C$ are equal, then $AB=c$, and the cutting off ``the lesser AC from the greater AB" cannot be carried out. On the other hand, if $AB=b$ and $b>c$, the triangle $ABC$  is not isosceles, and angles at $B$, $C$ are not equal. Throughout the proof, thus, the diagram \textit{changes} its metrical characteristic and cannot meet the assumptions of the proposition. Contrary to Euclid's claim, $D$ is a random point on $AB$, rather than introduced \textit{via} the following construction

  \begin{table}[!h]
 \centering
 \begin{tabular}{c}
$(B,c), AB$  \\ \hline
 $D$				 \\
 \end{tabular}
 \end{table}

  Now,  by SAS,  the equality $\triangle DBC= \triangle ACB$ holds, and Euclid concludes \textit{the lesser to the greater. The very notion is absurd}.

This time, the trichotomy law applies to triangles.
The contradiction
\[\triangle DBC=\triangle ACB\ \ \&\ \ \triangle DBC<\triangle ACB\]
occurs against the rule: \textit{For triangles, exactly one of the following conditions holds}
\[\triangle_1<\triangle_2,\ \ \triangle_1=\triangle_2,\ \ \ \triangle_1>\triangle_2.\]

   \begin{figure}
\begin{center}
\includegraphics[scale=0.7]{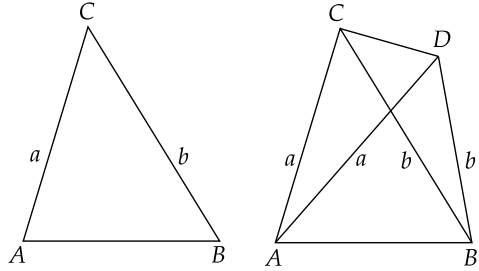}
\caption{\textit{Elements}, I.7 -- letters $a, b$ added} \label{figI7}
\end{center}
\end{figure}

\pv I.7 \textit{On the segment-line $AB$, two  segment lines   cannot meet at a different
point on the same side of \textit{AB}.}
 
The proof, atypically, includes no construction as the thesis explicitly states the impossibility of configuration depicted by the accompanying diagram. To get a contradiction, Euclid assumes there are two points $C, D$ such that $AC=a=AD$ and $BC=b=BD$ (see Fig. \ref{figI7}). 
Both triangles $\triangle ACD$ and  $\triangle BCD$ are isosceles and share the common base \textit{CD}. Angles at their bases  are equal, $\alpha=\alpha'$ and $\beta=\beta'$.

%\centerline{\includegraphics[scale=0.8]{I_7}}
\centerline{\includegraphics[scale=0.7]{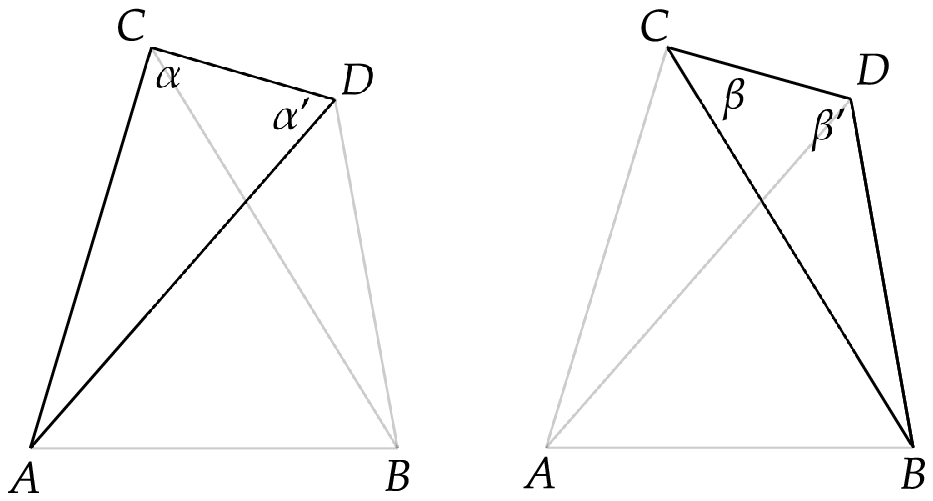}}
 
  Due to visual evidence,  at the vertex $C$, the inequality $\alpha>\beta$ holds, while at $D$,  $\beta'>\alpha '$.\footnote{\cite{ref_bmp} expounds the term \textit{visual evidence} in a bigger context.} Thus, $\beta'>\beta$ and, as stated earlier, $\beta'=\beta$. \textit{The very thing is impossible} -- clearly,  because exactly one of the conditions holds
\[\beta'<\beta,\ \ \beta'=\beta,\ \ \beta'>\beta.\]

That proof assumes the trichotomy rule for angles and transitivity of \textit{greater-than} relation. By modern standards, it is, thus, a total order.\footnote{Euclid applies the phrase ``is much greater than" when referring to the transitivity.}

In I.8, Euclid literally states the SSS criterion. Since the proof relies on a superposition of triangles,  we   propose the following paraphrase: 
 \textit{If two triangles share a common side and have other corresponding sides equal, then their corresponding angles will also be equal.}
In I.9--12, it is employed in that form as Euclid considers two equal triangles on both sides of the common side. %The construction of perpendicular plus SAS is enough to copy angles. 
Proof of that modification of I.8 effectively reduces to I.7.
%\footnote{Proof of I.7 gets complicated when point $D$ lies inside triangle $ABC$ (see Fig. \ref{figP}).} 
%(In the Hilbert system, transportation of angles enables proving SSS.)

\subsection{\textit{Greater-than} and \textit{Common Notions}}
Through \S\S\,10--11 of \cite{ref_RH}, Hartshorne seeks to prove Euclid's propositions I.1--34 
within the Hilbert system, except I.1 and I.23,  as they rely on the circle-circle axiom. 
He observes that ``Euclid's definitions, postulates, and common notions have been replaced by the undefined notions, definitions, and axioms"   in the Hilbert system. Commenting on Euclid's proofs of I.5--8, Hartshorne writes: ``Proposition I.5 and its proof is ok as they stand. [...] every step of Euclid's proof can be justified in a straightforward manner within the framework of a Hilbert plane. [...] Looking at I.6 [...] we have not defined the notion of inequality of triangles. However, a very slight change will give a satisfactory proof. [...] I.7 [...] needs some additional justification [...] which can be supplied from our axioms of betweenness [...]. 
For I.8, (SSS), we will need a new proof, since Euclid's method of superposition cannot be justified from our axioms" (\cite{ref_RH}, 97--99).

The above comparison between Euclid's and Hilbert's axiomatic approach simplifies rather than expounds. 
Euclid implicitly adopts \textit{greater-than} relation between line segments, angles, and triangles as primitive concepts; similarly to addition and subtraction (a lesser from the greater). In the previous section,  we have shown that he takes transitivity and the trichotomy law being self-evident. Further characteristics one can recover from his theory of magnitudes developed in Book V -- the only part of Euclid's geometry hardly discussed 
by Hartshorne (see \cite{ref_RH}, 166--167). Here is a brief account.

Euclidean proportion (for which we adopt the 17th-century symbol $::$) is a relation between two pairs of geometric figures (\textit{megethos}) of the same kind, triangles being of one kind, line segments of another kind, angles of yet another.
 Magnitudes of the same kind 
form an ordered additive semi-group \mbox{$\mathfrak{M}=(M,+,<)$} characterized by the five axioms given below (\cite{ref_BP}, \S\,3).

\begin{enumerate}\itemsep 0mm
\item[E1]$(\forall{a,b\in M})(\exists{n\in \N})(na>b)$.
\item[E2]$(\forall{a,b\in M})(\exists{c\in M})(a>b\Rightarrow a=b+c)$.
\item[E3] $(\forall{a,b,c\in M}) (a>b\Rightarrow{a+c>b+c})$.
\item[E4]$(\forall{a\in M})(\forall{n\in \N})(\exists{b\in M})(nb=a)$.
\item[E5]$(\forall{a,b,c\in M})(\exists{d\in M})(a:b::c:d), \ \ \mbox{where} \ \ na=\underbrace{a+a + ... + a}_{n - times}$.
\end{enumerate}

Clearly, E1--E3 provide extra characteristics of the \textit{greater-than} relation.

A modern interpretation of \textit{Common Notions} is simple:  CN 1 justifies the transitivity of congruence of line segments, triangles, and angles,
 CN 2 and 3 -- addition and subtraction in the following  form 
 \[a=a',\ b=b'\Rightarrow a+b=a'+b',\ \ a-a'=b-b'.\]
 
 The famous CN 5, \textit{Whole is greater than the part},  allows an interpretation by the formula $a+b>a$ (\cite{ref_bmp}, 73--76).

In the Hilbert system, the \textit{greater-than} relation  is defined through the concept of \textit{betweeness} and refers only to line segments and angles   (\cite{ref_RH}, 85, 95); similarly, addition of line segments and angles are introduced by definitions  (\cite{ref_RH}, 168, 93).  
Then counterparts of Euclid's axioms E2, E3, CN 1--3 are proved as theorems.
%Thus, Hilbert's attitude is rather \textit{genetic}, while Euclid's \textit{axiomatic} -- to employ \cite{ref_DH00}'s distinction portraying the 19th-century ways of introducing real numbers.\footnote{The status of \textit{greater-than} relation is even more contentious when viewed from the perspective of models; 
%in each Euclidean field, there is only one total order compatible with sums and products,  whereas there are Pythagorean fields that allow multiple orders that agree with sums and products.}  

\begin{figure}
\begin{center}
\includegraphics[scale=0.7]{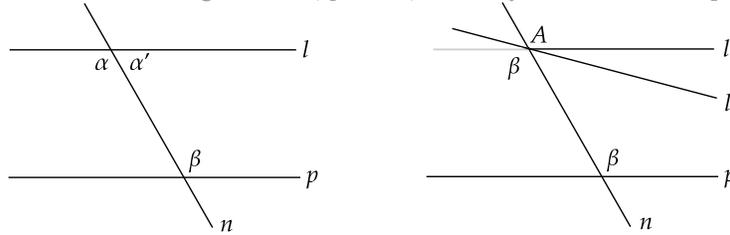}
\caption{\textit{Elements}, I.29 schematized (left). Hartshorne's version (right)} \label{figI29}
\end{center}
\end{figure}

Here is a sample argument based on inequalities and its Hilbert-style counterpart.
In I.29, Euclid proves  the thesis: \textit{When a line falls across parallel lines $l,\,p$,  equality of angles obtains} $\alpha =\beta$ (see Fig. \ref{figI29}). For, if they are not equal, one of the angles is greater, suppose  
$\alpha >\beta$. Then (implicitly by E3),%\footnote{In the enunciation of the proposition, instead of \textit{plus}, there is simply $AGH, BGH$, i.e., $\alpha, \alpha'$.}  
 \[\alpha >\beta \Rightarrow \alpha +\alpha'>\beta+\alpha'. \]
 
 Since $\alpha +\alpha'=\pi$, angles $\beta,\,\alpha'$ satisfy the requirement of the parallel axiom, i.e., 
        $\beta +\alpha' <\pi$ and straight lines $l,\,p$ meet, contrary to initial assumption.

{In contrast}, Hartshorne's proof of I.29 rests on the axiom stating there is exactly one line through the point $A$ parallel to $p$ (see Fig. \ref{figI29}). Then, if $\alpha\neq \beta$, a line $l'$ through $A$  making angle $\beta$ with $n$,  by I.27, is parallel to $p$ -- it contradicts the uniqueness of a parallel line through $A$.
 Euclid's proof, thus, implies an intersection point of $l$ and $p$, Hartshorne's -- a second parallel line to $p$.

\section{Existence}
\subsection{Existence \textit{via} Hilbert axioms} 

%Let us take a retrospective look at Euclid's constructions. 
Tables in sections \S\,1--2 expound grounds for introducing points, namely: a point is (1) an end of a given line segment, (2) a random point on a line segment, a straight line, a half-line, or a circle, (3) a circle-circle, circle-line, or line-line intersection point. 
Their existence  is covered by Hilbert axioms of \textit{Incidence}, \textit{Betweenness}, and Pasch, while circle-circle and circle-line intersection points require the circle-circle axiom.
%\footnote{In Greenberg's account, B-4 is the Pasch axiom.}
      Propositions I.7 and I.27 bring in other cases: (4) a random point on the plane, (5) a vertex of a triangle
that exists owning to the definition of parallel lines. 

Hilbert's characterization of a plane does not explain case I.7: the axiom on the existence of three non-colinear points 
allows to introduce triangle $ABC$ (see Fig. \ref{figI7}), yet, there are no grounds for point $D$ in the Hibert system. Similarly, Euclid assumes point $D$ while none of the above rules (1)--(3) guarantee its existence, and indeed, the very proposition does not include a construction part.  Commenting on I.7, Hartshorne points out the implicit argument on \textit{betweenness}, but does not report the suspicious  
\mbox{status of $D$  (see \cite{ref_RH}, 96, 99).} 

\begin{figure}
\begin{center}
\includegraphics[scale=0.7]{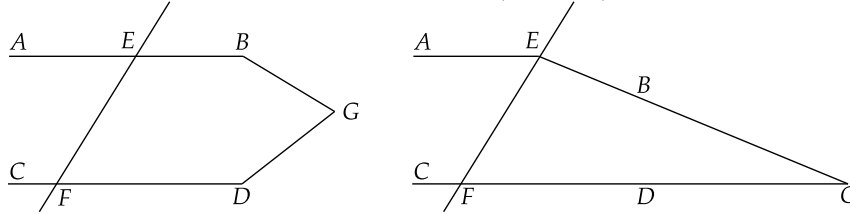}
\caption{\textit{Elements}, I.27 (left).  Triangle assumed in the proof (right).} \label{figI27}
\end{center}
\end{figure}

Ad I.27. Until proposition I.29, Euclid's arguments do not rely on the parallel postulate, yet, in I.27,   aiming to show $AB \parallel CD$, given that  $\angle AEF=\angle EFD$ (see Fig. \ref{figI27}), he invokes definition  of parallel lines: ``Parallel lines are straight-lines which, being in the
same plane, and being produced to infinity    in each direction,
meet with one another in neither" (I def. 23).

%\centerline{\includegraphics[scale=0.7]{I_27}}

The proof proceeds \textit{reductio ad absurdum} mode and starts with the claim: ``if not, being produced, AB and CD will certainly  meet together". Suppose, thus,  $AB$ and $CD$
 are not parallel and meet in $G$, then,  in triangle $EFG$, the external angle $\angle AEF$ is equal to the internal and opposite angle $\angle EFD$, but by I.16, $\angle AEF$ is also greater 
 than $\angle EFD$. Hence, $\angle AEF=\angle EFD$ and $\angle AEF>\angle EFD$. \textit{The very thing is impossible}.
 
The rationale for point $G$ lies in the definition of parallel lines rather than in
construction with straightedge and compass. $G$ is not, and cannot be determined as an intersection of lines $AB$ and $CD$: straight lines $AB$ and $CD$ do not meet in either absolute or Euclidean geometry.

Summing up, 
points $D$, in I.7, and $G$, in I.27, constitute impossible figures that can not be constructed with a straightedge and compass. Similar arguments 
apply to the figure I.39 and those in Book III accompanying propositions 2,\,10,\,13,\,16,\,23,\,24. Generally, such impossible figures reveal Euclid's struggle with foundational problems: I.39  aims to link the theory of equal figures and parallel lines, those in Book III -- making the foundations of trigonometry.

\subsection{Existence and \textit{co-exact attributes}}

\begin{figure}
\begin{center}
\includegraphics[scale=0.7]{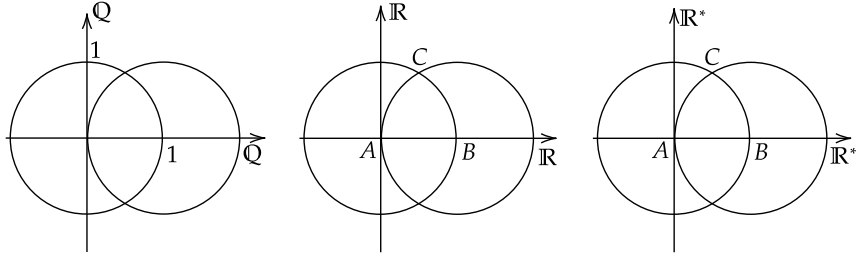}
\caption{Circles on non-Euclidean and Euclidean planes} \label{fig3o}
\end{center}
\end{figure}

Euclid proposition I.1 is a model example for the proponents of diagrammatic thinking.
The enunciation of the proposition introduces point $C$ depicted on the accompanying diagram through the following phrase  ``the point C, where the circles cut one another" (see Fig. \ref{fig3o}). 
It is not the case one has to read off from the diagram that involved circles intersect. The circle-circle axiom is not explicitly assumed, but it is not the only tacit supposition of the \textit{Elements}.  
Modern geometry studies how it relates to other axioms and how it affects the characterization of a plane. The key result in that respect states: if $\F\times\F$ is a  Cartesian plane  over an ordered field $(\F,+,\cdot,0,1,<)$, the field
is closed under the square operation iff the circle-circle axiom is satisfied  on $\F\times\F$ (\cite{ref_RH}, 144).

 Let us consider now Manders' claim on co-exact attribute as it explicitly refers to I.1. It reads:
``Co-exact attributes
are those conditions which are unaffected by some range of every continuous
variation of a specified diagram [...] or the existence of intersection points such as
those required in Euclid I.1 (which is unaffected no matter how the circles
are to some extent deformed)" (\cite{ref_KM}, 92).

The key concepts here are \textit{deformation} and \textit{continuous variation}. Manders suggest reference to topological \textit{deformations}, yet does not provide any details. Taken literally, they require specification of the plane on which a diagram lays.
{Instead of any \textit{deformation}}, we  consider diagram I.1 on three planes: $\Q\times \Q$, $\R\times \R$, and $\R^*\times \R^*$  (see Fig. \ref{fig3o}).\footnote{$\R^*$ is the set of hyperreals (see \S\,6 below); circles with radius $1$ have centers at $(0,0)$ and $(0,1)$.} The second and third are models of the Euclidean planes as fields of real and hyperreal numbers are closed under the square root operation, while the first is not a  Euclidean plane. Due to calculations, we get to know that on the plane $\Q\times \Q$ the intersection point of circles does not exist. However, from the cognitive perspective, one and the same diagram represents circles in these three various mathematical contexts.\footnote{As for $\Q\times \Q$, we get sure that circles look like continuous objects due to a technical result such as  \cite{ref_LT}.} Thus, with no deformation 
of the diagram, 
but  by switching from one mathematical context to another we get or not that the intersection point of the circles exists.

\section{Inequalities and \textit{co-exact} attributes}

In his analysis of proposition I.6, Manders interprets Euclid's argument $\triangle DBC$ is smaller than
$\triangle ABC$ in terms of part-whole relation read-off from the diagram (\cite{ref_KM}, 109--110). On another occasion, he explicitly writes: ``A strict inequality may be reduced [...] to a proper-part
relationship in the diagram" (\cite{ref_KM}, 112).

\begin{figure}[h!]
\begin{center}
\includegraphics[scale=0.7]{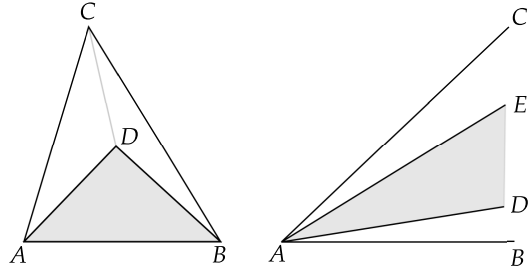}
\caption{Non-Eculidean parts} \label{figP}
\end{center}
\end{figure}

In \cite{ref_bmp}, we provide a detailed analysis of Common Notion 5  showing    it does not reduce to part-whole relation. In short, logical analysis gives the following formula for CN5: $a+b>a$, where $a+b$ stands for \textit{the whole} and $b$ -- for \textit{its part}, and $a+b,\,a,\,b$ have to be of the same kind (triangles or angles). Given that interpretation, triangle $ADB$ is not a Euclidean part of triangle $ACB$; similarly, the gray angle is not a Euclidean part of the angle $CAB$, \mbox{as represented in Fig. \ref{figP}.} 

Euclid's I.32 sets another challenge for the diagrammatic philosophy.  Fereir{\'o}s reads it out as follows:   ``Most of the proof steps involve exact information, dealing with equalities of angles, and they depend on previous theorems (I.31, I.29) and common notions 2. Only two steps involve co-exact information, bringing in attributions read from the diagram", namely,  that $\angle ACD$ covers $\angle AGE$ and $\angle EGD$ 
``is co-exact because it has to do with part-whole relation and it is not affected by deformation" (\cite{ref_JF}, 135--136).

\begin{figure}
\begin{center}
\includegraphics[scale=0.7]{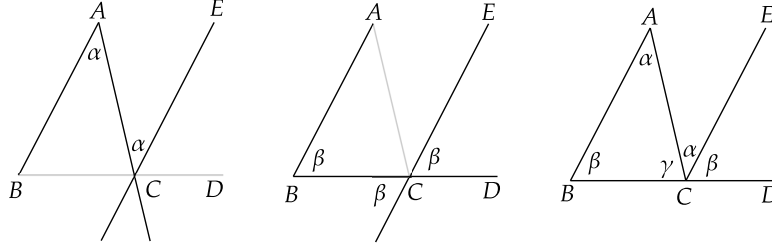}
\caption{Proof of I.32 schematized.} \label{figI32c}
\end{center}
\end{figure}

Here is our interpretation (see Fig. \ref{figI32c}). Euclid transports angle $\alpha$ to point $C$, and draws $CE$, which,  by I.27, is parallel to $AB$. Hence, by I.29, $\angle ECD=\beta$, and angles at $C$ sum up to ``two right angles", $\gamma+\alpha+\beta=\pi$.

 Now, $CE$ lies inside angle $ACD$ because, by I.17,
$\gamma+\alpha<\pi$. It is not the case, thus, that the only way to get that information is to read it from the diagram. 
\section{Existence meets Inequalities}

\subsection{Semi-Euclidean plane}

In this section, we present a model o semi-Euclidean plane, i.e., a plane in which angles in a triangle sum up to $\pi$ yet the parallel postulate fails. 
\cite{ref_RH}, p. 311, introduces that term, but the very idea originates in Max Dehn's  1900 \cite{ref_MD}, \S\,9, which built such a model owning to a non-Archimedean Pythagorean field. Dehn explored a non-Euclidean field introduced already in Hilbert \cite{ref_DH}, \S\,12.\footnote{See also \cite{ref_RH}, \S\,18. Example 18.4.3 expounds Dehn's model.} We employ the Euclidean field of hyperreal numbers. On the Cartesian plane over hyperreals, the circle-circle and circle-line intersection axioms are satisfied, meaning one can mirror Euclid's straightedge and compass constructions. To elaborate, let us start with the introduction of the hyperreal numbers.

%\subsection{Hyperreals}

An ordered field $(\F,+,\cdot,0,1,<)$ is a  commutative field together with a total order that is compatible with sums and products.
%\[x<y \Rightarrow x+z<y+z,\ \ \ x<y,\, 0<z\Rightarrow xz<yz.\]
 In such a field, one can define the following subsets of $\F$:
 \begin{enumerate}\itemsep 0mm
 \item[] $\mathbb L= \{x\in \F: (\exists n\in \N)(|x|<n)\}$,
 \item[] $\Psi= \{x\in \F: (\forall n\in \N)(|x|>n)\}$,
 \item[] $\Omega = \{x\in \F: (\forall n\in \N)(|x|<\tfrac 1n)\}$.
 \end{enumerate}

They are called limited, infinite, and infinitely small numbers, respectively.
 Here are some relationships helpful to pursue our arguments.
\begin{enumerate}\itemsep 0mm
\item[] $(\forall x, y\in\Omega)( x+y\in\Omega, x y\in\Omega)$,
\item[]  $(\forall x\in\Omega)( \forall y\in\mathbb L)( x y\in\Omega)$,
%\item[] $(\forall x)(x\in \mathbb A\Rightarrow  x^{-1}\in\mathbb A)$,
\item[] $(\forall x\neq 0)(x\in\Omega  \Leftrightarrow \ x^{-1}\in\Psi)$.
\end{enumerate}

To clarify our account, let us observe the following equality $\Omega=\{0\}$ is a version of the well-known    Archimedean axiom.
Since real numbers form the \textit{biggest} Archimedean field, every field extension of $(\R,+,\cdot,0,1,<)$
includes positive infinitesimals. 

\begin{figure}
\begin{center}
\includegraphics[scale=0.8]{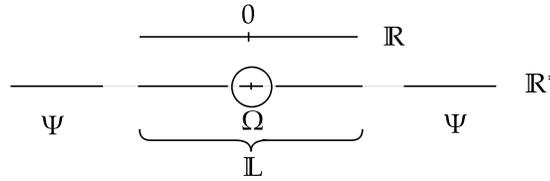}
\caption{The line of real numbers and its extension to hyperreals} \label{fig5}
\end{center}
\end{figure}

Let $\mathcal U$ be a non-principal ultrafilter on $\N$.
The set of hyperreals is defined as a reduced product $\R^*=\R^{\N}/{\mathcal U}$. Sums, products, and the ordered are introduced pointwise.
A reader can take for granted that 
the field of hyperreals $(\R^*,+,\cdot,0,1,<)$
extends real numbers, hence, includes infinitesimals and infinite numbers; moreover, it is closed under the square root operation (see \cite{ref_PB16}, \cite{ref_PB01}). 
Fig. \ref{fig5} represents in a schematized way a relationship between $\R$ and $\R^*$, as well as between $\mathbb L$, $\Psi$, and $\Omega$. 

Due to the proposition 16.2 (\cite{ref_RH}, 144), the Cartesian plane over the field of hyperreals is a model of Euclidean plane, with straight lines and circles given by equations $ax+by+c=0$, $(x-a)^2+(y-b)^2=r^2$, 
where $a,b,c,r\in \R^*$; angles between straight lines are defined as in the Cartesian plane over the field of real numbers. Specifically, on the plane $\R^*\times\R^*$, angles in triangles  sum up to $\pi$. Parallel lines are of the form $y=mx+b$ and $y=mx+c$, while a perpendicular to the line $y=mx+b$ is given by the formula 
$y=-\frac{1}{m}x+d$.

Now, take us a subspace $\mathbb L\times\mathbb L$ of the plane  $\R^*\times\R^*$. On that plane, circles are defined by analogous formula, namely 
$(x-a)^2+(y-b)^2=r^2$, where $a,b,r\in \mathbb L$, while  every line in $\mathbb L\times\mathbb L$ 
is of the form $l\cap \mathbb L\times\mathbb L$, where $l$ is a line in $\R^*\times\R^*$. %That modification finds the following substantiation: 
Since we want plane $\mathbb L\times\mathbb L$ include  lines such as $y_1=\varepsilon x$, where  
$\varepsilon\in\Omega$,  it has also to include the perpendicular $y_2=\frac{-1}{\varepsilon} x$, but 
$\frac{-1}{\varepsilon}\notin\mathbb L$. Formula 
$l\cap \mathbb L\times\mathbb L$, where $l=ax+by+c$ and $a, b, c \in \R^*$   guarantees the existence of the straight line $y_2$ in  $\mathbb L\times\mathbb L$. 
Finally, the interpretation of angle is the same as in the model $\R^*\times \R^*$.
\begin{figure}
\begin{center}
\includegraphics[scale=0.8]{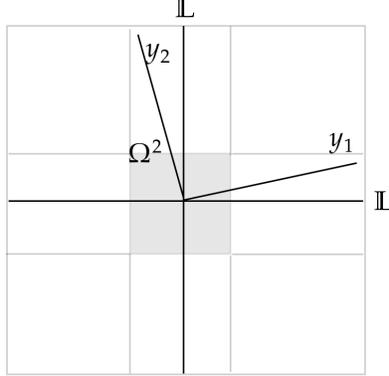}
\caption{Perpendicular lines with infinitesimal and infinitely large slopes} \label{fig6}
\end{center}
\end{figure}

Explicit checking shows that the model characterized above satisfies all Hilbert axioms of non-Archimdean  plane geometry plus the circle-circle and line-circle axioms, except parallel axiom; 
the more general theorem concerning Hilbert planes also justifies our model, namely \cite{ref_RH}, p. 425, theorem, 43.7 (a).

\begin{figure}
\begin{center}
\includegraphics[scale=0.7]{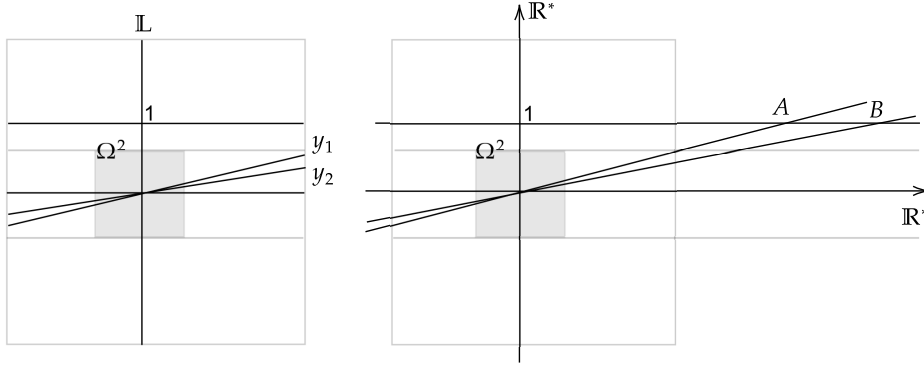}
\caption{Non-Euclidean plane $\mathbb L\times \mathbb L$ (left) vs. Euclidean plane $\R^*\times \R^*$ (right)} \label{fig7}
\end{center}
\end{figure}

With regard to parallel lines, let us consider the horizontal line $y=1$ and two specific lines through $(0,0)$,
$y_1=\varepsilon x, y_2=\delta x$, where $\varepsilon, \delta\in\Omega$. (see Fig. \ref{fig7}).
Since
$\Omega\mathbb L\subset \Omega$, the following inclusions hold $y_1, y_2\subset \mathbb L\times \Omega$. In other words, values of maps $y_1(x), y_2(x)$ are infinitesimals, given that $x\in\mathbb L$.  The same obtains for any line of the form $y=\mu x$, with $\mu\in\Omega$. Since there are infinitely many infinitesimals, 
there are infinitely many lines through $(0,0)$ not intersecting the horizontal line $y=1$.

\begin{figure}
\begin{center}
\includegraphics[scale=0.7]{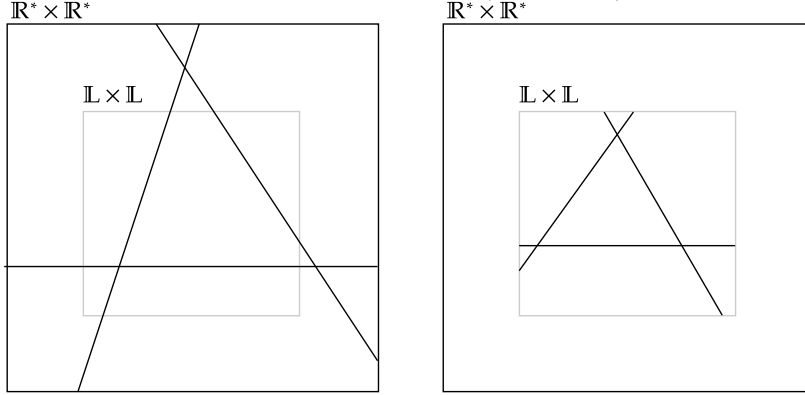}
\caption{Triangles in Euclidean plane $\R^*\times \R^*$ and its subspace $\mathbb L\times \mathbb L$} \label{figT}
\end{center}
\end{figure}

Since every triangle in $\mathbb L\times \mathbb L$ is a triangle in $\R^*\times \R^*$, it follows that angles in a triangle  on the plane $\mathbb L\times \mathbb L$  sum up to $\pi$ (see Fig. \ref{figT}).

\subsection{Exact attributes on the semi-Euclidean plane}

Let us go back to Manders' claim ``Exact attributes are those which, for at least
some continuous variation of the diagram, obtain only in isolated cases;
paradigmatically, straightness of lines or equality of angles (neither of which
survive any except exceptional types of deformation, no matter how small)". 

In Fig. \ref{fig7}, points  $A=(\frac 1{\varepsilon},1)$, $B=(\frac 1{\delta},1)$  are intersections  of $y_1, y_2$ with the horizontal line $y=1$ on the plane $\R^*\times\R^*$. Vertical line $x=0$ falls on $y=1$ and $y_1$ \textit{making internal angles}   less than $\pi$. Switching from $\R^*\times\R^*$ to   $\mathbb L\times\mathbb L$, we do not modify  \textit{straightness of lines or equality of angles} as far as it regards the diagram on 
$\mathbb L\times\mathbb L$, however,
 there is nothing in the diagram which informs us whether they intersect or not. Knowing plane characteristics, we can infer they intersect, given that being on $\R^*\times\R^*$. {The clue is the diagrammatic philosophy does not impose any restrictions on Euclid's clause ``being produced to infinity".  
As focused on a diagram, it does not have means to decide whether lines meet beyond a diagram. 

Possibly, the claim on \textit{straightness of lines or equality of angles} is designed to exclude hyperbolic representations of lines and angles. Indeed, in the Poincare model, straight lines are circles and angles are determined between tangents; in the Klein model, straight lines are straight, while angles are retrieved from the Poincare model. Both models then change the diagrammatic representation of straight lines and angles.  In the plane $\mathbb L\times\mathbb L$, straight-lines are usual straight-lines, angles are usual angles; moreover, triangles are Euclidean, there are also squares and usual perpendicular lines.  Yet, there are no instruments in the diagrammatic tool-box to reveal a non-Euclidean character of that plane.  
It is because diagrammatic perspective is local, while the phenomenon of parallelism requires a global perspective.}

%``Exact attributes (indeed, by definition) are unstable under perturbation of a
%diagram" (\cite{ref_KM}, 93).

%\section{Summary}
%Euclidean diagram, likewise a modern mathematical symbol, be it an integral $\int_a^bf(x)dx$ or $e^{ix}$,
%does not speak for itself -- it requires context. Manders suggests studying \textit{continuous variation of a diagram} will do; we sought to argue that techniques of modern mathematics enable one to dig deeper into  the \textit{Elements}.

%
%
%\section{Acknowledgments}
%The first author is supported by the National Science Centre, Poland grant
%2018/31/B/HS1/03896. The second author is supported by the National Science Centre, Poland grant 2020/37/N/HS1/01989.

\end{document}